\documentclass[12pt]{article}
\usepackage{amssymb, amsmath}

\newtheorem{theorem}{Theorem}[section]
\newtheorem{lemma}{Lemma}[section]

\newtheorem{cor}{Corollary}[section]

\newcommand{\n}{\nonumber}

\newcommand{\si}{\sigma_R }

\newcommand{\s}{\sigma}

\renewcommand{\o}{\omega}
\newcommand{\q}{\mathcal{Q}}

\newcommand{\bb}{\begin{equation}}
\newcommand{\ee}{\end{equation}}
\newcommand{\bq}{\begin{eqnarray}}
\newcommand{\eq}{\end{eqnarray}}
\newcommand{\bqn}{\begin{eqnarray*}}
\newcommand{\eqn}{\end{eqnarray*}}

\begin{document}
\title{Localized energy equalities  for the Navier-Stokes and the Euler equations}
\author{ Dongho Chae \\
 Department of Mathematics\\
 Chung-Ang
University\\
 Seoul 156-756, Korea\\
 e-mail: {\it dchae@cau.ac.kr }}
 \date{}
\maketitle
\begin{abstract}
Let  $(v,p)$ be a smooth solution pair of the velocity and the
pressure for the Navier-Stokes(Euler) equations on $\Bbb R^N\times
(0, T)$, $N\geq 3$.  We set the Bernoulli function $Q=\frac12 |v|^2
+p$. Under suitable decay conditions at infinity for $(v,p)$ we
prove that for almost all $\alpha(t)$ and $\beta(t)$  defined on
$(0, T)$ there holds
 \bqn
&&\int_{\{\alpha(t)< Q(x,t)<\beta(t)\}}
\left(\frac12\frac{\partial}{\partial t} |v|^2+\nu |\o|^2
 \right) dx
 =\nu \int_{\{ Q(x,t)= \beta(t)\}} |\nabla Q|dS\\
 &&\hspace{1.7in}-\nu \int_{\{ Q(x,t)= \alpha(t)\}} |\nabla Q|dS,
 \eqn
where $\o=$ curl $v$ is the vorticity.
 This shows  that, in each region squeezed between two levels of the Bernoulli
function, besides the energy dissipation due to the enstrophy, the
energy flows into the region through the level hypersurface having
the higher level, and the energy flows out of the region through the
level hypersurface  with the lower level.  Passing
$\alpha(t)\downarrow\inf_{x\in \Bbb R^N} Q(x,t)$ and $
\beta(t)\uparrow\sup_{x\in \Bbb R^N} Q(x,t)$, we recover the
well-known energy equality, $
   \frac12\frac{d}{dt} \int_{\Bbb R^N}  |v|^2=-\nu\int_{\Bbb R^N} |\o|^2 dx.$
  A weaker version of the above equality under the weaker decay assumption of the solution at spatial infinity
 is also derived. The stationary version of the equality implies the previous
 Liouville type results on the Navier-Stokes equations. 
\end{abstract}
\section{Introduction}
 \setcounter{equation}{0}
We are concerned  here  on the  incompressible Navier-Stokes(Euler)
equations on $\Bbb R^N$, $N\geq 2$.
  $$
(NS,E)\left\{ \aligned & v_t+ ( v\cdot \nabla )  v +\nabla p=\nu\Delta v,\\
&\mathrm{div} \, v=0,\\
 &v(x,0)=v_0 (x),
 \endaligned
 \right.
 $$
 where $v=v(x,t)=(v_1(x,t),\cdots , v_N(x,t))$ is the velocity, and  $p=p(x,t)$ is the pressure.
 We assume the viscosity satisfies $\nu \geq 0$. In the case $\nu>0$  the system (NS,E) becomes
 the Navier-Stokes equations, while
 for $\nu =0$ the system (NS,E) is the Euler equations.
By the  system (NS,E) we represent both of the cases of the
Navier-Stokes and the Euler equations. Below we denote the vorticity
of the vector field $v$ in $\Bbb R^N$ defined by
  $\o =\{ \partial_j v_k-\partial_k v_j\}_{j,k=1; j>k}^N$, the
  magnitude  of which is given by
  $$ |\o|=\sqrt{\frac12 \sum_{j,k=1}^N(\partial_j v_k-\partial_k v_j)^2}.$$
  As is well-known in most of the fluid mechanics text books, the
  smooth solution $v$ of the system (NS,E) with sufficiently fast decays at spatial infinity
  satisfies the energy equality,
  \bb\label{energy}
  \frac12 \frac{d}{dt} \int_{\Bbb R^N} |v|^2 dx=-\nu \int_{\Bbb R^N}
  |\o |^2\, dx.
  \ee
Our aim in this paper is to localize the domain of the integration
in (\ref{energy}). The domains are characterized by their
boundaries, which are the level hyper-surfaces of the Bernoulli
function, \bb\label{ber}
 Q(x,t):=\frac12 |v(x,t)|^2 +p(x,t).
\ee
 Our
localized integral equalities refine the classical equality
(\ref{energy}), in the sense that under suitable integrability
conditions for the solutions and particular choice of the levels of
the Bernoulli  function  we recover (\ref{energy}).
  The first theorem below concerns on these localized energy equalities
  under milder conditions on the asymptotic behavior for the solution $(v,p)$ at spatial infinity.
\begin{theorem}
Let $N\geq 2$, and $(v, p)$ be a smooth solution of (NS,E) on $\Bbb
R^N\times (0, T)$. Let $Q(x,t)$ be defined as in (\ref{ber}).
 Suppose there exists $Q_0=Q_0
(t)$ such that
 \bb
 \lim_{|x|\to \infty} Q(x,t)=Q_0(t)
 \ee
 uniformly for each $t\in (0, T)$. Then, for
almost every $\alpha(t),\beta(t)$ defined on $(0, T)$ such that
either $\beta (t) >\alpha(t)>0$ or $ 0>\beta (t)
 >\alpha(t)$
we have the following equality.
 \bq\label{th1}
\lefteqn{\int_{\{ \alpha(t)<Q(x,t)- Q_0(t)<\beta(t)\}} \left(\frac12
\frac{\partial}{\partial t} |v|^2 +\nu|\o|^2
 \right) dx}\n \\
 && \quad=\nu \int_{\{ Q(x,t)- Q_0(t)=\beta(t)\}} |\nabla Q
 |dS-\nu \int_{\{ Q(x,t)- Q_0(t)=\alpha(t)\}} |\nabla Q
 |dS.\n \\
 \eq
  \end{theorem}
 {\em Remark 1.1 }   Physically the equality (\ref{th1}) shows new energy balance in each region
 squeezed between two levels of the Bernoulli
function.  Besides the well-known energy dissipation due to the
enstrophy the energy flows into the region through the level
hypersurface having the higher level, and the energy flows out of
the region through the level hypersurface  with the lower level.
This tendency of flow of the energy into lower values of the
Bernoulli function is proportional to the viscosity and the
magnitude of the gradient of the Bernoulli function,
 and there exists no such phenomena in the inviscid case of the Euler equations, where
  energy is just conserved on the average in each strip between
the two levels of the Bernoulli function.
 \\
 \ \\
 {\em Remark 1.2 } If we approach $\beta(t)\uparrow \sup_{x\in \Bbb R^N}
 Q(x,t)-Q_0(t)$, then,  since $\int_{\{ Q(x,t)- Q_0(t)=\beta(t)\}} |\nabla Q
 |dS\to 0$,  we obtain from (\ref{th1}) that
 \bb\label{th1aa}
\int_{\{ \alpha(t)+Q_0(t)<Q(x,t) \}} \left(\frac12
\frac{\partial}{\partial t} |v|^2 +\nu|\o|^2
 \right) dx
=-\nu \int_{\{ Q(x,t)=\alpha(t)+Q_0(t)\}} |\nabla Q
 |dS \leq 0.
 \ee
Similarly, if we approach $\alpha(t)\downarrow\inf_{x\in \Bbb R^N}
 Q(x,t)-Q_0(t)$, then we obtain from (\ref{th1}) that
 \bb\label{th1bb}
\int_{\{ Q(x,t)<\beta(t)+Q_0(t)\}} \left(\frac12
\frac{\partial}{\partial t} |v|^2 +\nu|\o|^2
 \right) dx
=\nu \int_{\{ Q(x,t)=\beta(t)+Q_0(t)\}} |\nabla Q
 |dS \geq 0.
 \ee
 {\em Remark 1.3 } If we let $\alpha(t)$ and $\beta(t)$ approach respectively to the critical values of
 $Q(x,t)-Q_0(t)$ with $\alpha(t)\beta(t) >0$,
 then the righthand side of (\ref{th1}) vanishes, and the
 energy balance holds in the region of the set $\{ x\in
 \Bbb R^N\, |\, \alpha(t)<Q(x,t)- Q_0(t)<\beta(t)\}$.\\
\\
The following theorem derives stronger version of equalities than
(\ref{th1}) in the sense that
 there is no sign condition on the functions $\alpha(t), \beta(t)$, but under the stronger assumptions
 on the behaviors of the solutions at spatial infinity represented by the integrability of the solutions.
\begin{theorem} Let $N\geq 3$, and  $(v, p)$ be a smooth solution of (NS,E) on $\Bbb R^N\times (0,
T)$.
 Suppose $(v,p)$
satisfies the following conditions: there exists $p_0(t)$ defined on
$(0, T)$ such that
 \bb\label{hypo1}
 \lim_{|x|\to \infty} (|v(x,t)|+|p(x,t)-p_0(t)|)=0,
 \ee
 uniformly for each $t\in (0, T)$,
\bb\label{hypo}
 v\in L^\infty (0, T; L^{\frac{3N}{N-1}} (\Bbb R^N)\cap L^{\frac{2N}{N-2}} (\Bbb R^N)),
 \ee
 and
 \bb\label{hypo2}
\frac{\partial}{\partial t} |v|^2 \in L^\infty (0, T; L^1 (\Bbb
R^N)).
 \ee
Then, for almost all real valued functions $\alpha(t)<\beta(t)$
defined on $(0, T)$ we have the following equalities.
 \bq\label{th2}
\lefteqn{\int_{\{ \alpha(t)+p_0(t)<Q(x,t)<\beta(t)+p_0(t)\}}
\left(\frac12 \frac{\partial}{\partial t} |v|^2 +\nu|\o|^2
 \right) dx}\n \\
 && \quad=\nu \int_{\{ Q(x,t)=\beta(t)+p_0(t)\}} |\nabla Q
 |dS-\nu \int_{\{ Q(x,t)=\alpha(t)+p_0(t)\}} |\nabla Q
 |dS.\n \\
 \eq
  \end{theorem}
\noindent{\em Remark 1.4 }  The same remark holds as Remark 1.1 for
the physical interpretation of (\ref{th2}).  If we approach
    $$\alpha(t)\downarrow\inf_{x\in \Bbb R^N} Q(x,t)-p_0(t),\quad \beta(t)\uparrow\sup_{x\in \Bbb R^N} Q(x,t)-p_0(t), $$
    in particular, then from (\ref{th2}) we obtain
    \bb\label{ee}
    \int_{\Bbb R^N} \left(\frac12\frac{\partial}{\partial t} |v|^2 +\nu |\o|^2
 \right) dx=0,
 \ee
 which is equivalent to  the classical energy equality (\ref{energy}) for the Navier-Stokes equations.\\

 For the stationary solutions of the Navier-Stokes equations with $v(x,t)=v(x), p(x,t)=p(x),$ and $\o(x,t)=\o(x)$,
  as
 an immediate corollary of Theorem 1.2 we obtain the following
 result, which is previously obtained by Galdi(\cite{gal}).
 \begin{cor} Let $N\geq 3$.
 Suppose $(v,p)$ is a stationary smooth solution of the Navier-Stokes equations
 satisfying: there exists a constant $p_0$ such that
 \bb\label{sta1}
 \lim_{|x|\to \infty} (|v(x)|+|p(x)-p_0|)=0,
 \ee
 uniformly, and
\bb\label{sta2}
 v\in \left\{ \aligned & L^{\frac{9}{2}} (\Bbb R^3)\cap L^6(\Bbb R^3) \quad \mbox{if $ N=3$},\\
 & L^{\frac{2N}{N-2}}(\Bbb R^N) \quad \mbox{if $ N\geq 4 $.}\endaligned
 \right.
 \ee
 Then,  $v=0$.
 \end{cor}
Indeed, we first note that  $\frac{3N}{N-1} \geq \frac{2N}{N-2}$ for
$N\geq 4$. Thus,  if $N\ge 4 $,
 then $ L^{\frac{3N}{N-1}}
  (\Bbb R^N) \cap L^\infty(\Bbb R^N)\hookrightarrow L^{\frac{2N}{N-2}} (\Bbb R^N)\cap L^\infty (\Bbb R^N)$ by the standard $L^p-$interpolation, and
the fact $v\in L^\infty (\Bbb R^N)$ is guaranteed by the hypothesis
(\ref{sta1}).
 Therefore, the stationary version of (\ref{ee}) together with the assumptions (\ref{sta1})
 and (\ref{sta2}) imply $\o=0$. Combining this with div $v=0$, we have $v=\nabla h$ for some harmonic function $h$.
 Therefore the condition (\ref{sta1}) implies $v=0$. \\
\ \\
{\em Remark 1.5 } For $N=3$, if we assume (\ref{sta1}), and replace
the assumption (\ref{sta2}) by $\int_{\Bbb R^3} |\o|^2 dx <\infty$,
it is still an open question whether $v=0$ or not. In this case we
know by the maximum principle from (\ref{le}) implies that $Q(x)\leq
p_0 $ for all $x\in \Bbb R^3$. Therefore, by choosing
$\alpha\downarrow \inf_{x\in \Bbb R^3} Q(x)-p_0$, and $\beta
\uparrow 0$  in the stationary version of Theorem 1.1,  we obtain
 \bb
 \int_{\Bbb R^3} |\o|^2 dx= \lim_{\beta \uparrow 0} \int_{ \{Q
 (x)=p_0 +\beta\}} |\nabla Q|dS.
 \ee
Therefore, we find that the desired Liouville type theorem holds if
$$ \lim_{\beta \uparrow 0} \int_{ \{Q
 (x)=p_0 +\beta\}} |\nabla Q|dS < \int_{\Bbb R^3} |\o|^2 dx. $$
\ \\
In the case of smooth solutions to (NS,E) in a periodic domain
(\ref{sta1}) and (\ref{sta2}) are not necessary in view of the proof
of Theorem 1.2. Therefore, as another corollary of Theorem 1.2 and
its proof we obtain the following:
\begin{cor} Let $N\geq 2$.
 Suppose $(v,p)$ is a smooth solution of (NS,E)  on $\Bbb T^N
 \times (0, T)$, where $\Bbb T^N=\Bbb R^N/ \Bbb Z^N$ is the periodic
 domain in $\Bbb R^N$. We normalize the pressure by assuming $\int_{\Bbb T^N} p \,dx=0$.
 Define the Bernoulli function $Q$ as in (\ref{ber}). Then,
 for all $t\in (0, T)$ and for almost all real valued functions  $\alpha(t), \beta(t)$ we have
 Then, for all real valued functions $\alpha(t)<\beta(t)$ defined on
$(0, T)$ we have the following equalities.
 \bq\label{th2}
\lefteqn{\int_{\{ \alpha(t)<Q(x,t)<\beta(t)\}} \left(\frac12
\frac{\partial}{\partial t} |v|^2 +\nu|\o|^2
 \right) dx}\n \\
 && \quad=\nu \int_{\{ Q(x,t))=\alpha(t)\}} |\nabla Q
 |dS-\nu \int_{\{ Q(x,t)=\beta(t)\}} |\nabla Q
 |dS.\n \\
 \eq
 \end{cor}
\noindent{\em Remark 1.6 }  Similarly to Remark 1.4 the usual energy
equality in $\Bbb T^N$ follows immediately by approaching
    $$\alpha(t)\downarrow\inf_{x\in \Bbb T^N} Q(x,t),\quad\beta(t)\uparrow\sup_{x\in \Bbb T^N} Q(x,t).
    $$
    \ \\
   We emphasize that we consider the energy equalities only for smooth solutions of
    the Navier-Stokes and the Euler equations. For the studies of the energy for the weak solutions of
    the Euler equations we refer to the results, among others, by De
    Lelilis-Sz\'{e}kelyhidi(\cite{del}), Shnirelman(\cite{shn}) and the references therein.
\section{Proof of the Main Theorems }
\setcounter{equation}{0}
 Although the lemma below is well-known previously for the
 stationary solutions of the Navier-Stokes equations, we present its proof
 for reader's convenience.
\begin{lemma}
If $(v,p)$ solves (NS,E), then $Q=p+\frac12 |v|^2$ satisfies
 \bb\label{le}
 \frac12\frac{\partial}{\partial t} |v|^2 + \nu |\o|^2=\nu\Delta Q-v\cdot \nabla Q.
 \ee
\end{lemma}
{\bf Proof } Multiplying (NS,E) by $v$, we obtain
 \bb\label{le2}
 \frac12\frac{\partial}{\partial t} |v|^2+ v\cdot \nabla \left(p+\frac12 |v|^2\right)=\nu v\cdot \Delta v=\nu \Delta
  \left(\frac12|v|^2 \right)-\nu \sum_{j,k=1}^N ( \partial_j v_k )^2.
  \ee
  We compute
  \bq\label{le3}
\sum_{j,k=1}^N ( \partial_j v_k )^2&=& \sum_{j,k=1}^N  \partial_j
v_k (\partial_jv_k -\partial_k v_j )+\sum_{j,k=1}^N \partial_j v_k
\partial_k v_j \n \\
&=& \frac 12\sum_{j,k=1}^N(\partial_jv_k -\partial_k v_j )^2 -
\Delta p=|\o|^2 -\Delta p,
 \eq
where we used the well-known formula $\Delta p=-\sum_{j,k=1}^N
\partial_j v_k
\partial_k v_j$, which is obtained from (NS,E) by taking div$(\cdot)$.
Plugging (\ref{le3}) into (\ref{le2}), we have (\ref{le}).
$\square$\\
\ \\
\noindent{\bf Proof of Theorem 1.1 } We use argument similar to the
one used in \cite{cha}.   We define
$\mathcal{Q}(x,t):=Q(x,t)-Q_0(t)$. Then,
  \bb\label{13}
\lim_{|x|\to \infty} \mathcal{Q}(x,t)=0\quad\forall t\in (0, T)
  \ee
  uniformly.
From (\ref{le}) we have
 \bb\label{prin}
 \nu \Delta \mathcal{Q}-v\cdot \nabla \mathcal{Q}=\frac12
 \frac{\partial}{\partial t}|v|^2+\nu|\o|^2.
 \ee
We assume below
$$\sup_{x\in \Bbb R^N} \mathcal{Q}(x,t)\geq \beta
(t)
>\alpha(t)>0.
$$
 The other case, $0>\beta (t) >\alpha(t)\geq \inf_{x\in \Bbb R^N} \mathcal{Q} (x,t)$, is similar. Let us set
 $\Bbb D=\{ x\in \Bbb R^N\, |\,  \alpha(t) <\q (x,t)<\beta(t) \}.
 $
 Thanks to the Sard theorem and the implicit function theorem
$\partial \Bbb D$ consists of smooth level surfaces in $\Bbb R^N$
except the values of $\alpha(t)$ and $\beta(t)$, having the zero
Lebesgue measure, which corresponds to the critical values of the
function $z=\mathcal{Q}(x,t)$. It is understood that our values of
$\alpha(t), \beta(t)$ below avoid these exceptional ones.
 We integrate (\ref{prin}) on $\Bbb D$,
  and use the divergence theorem to obtain
 \bq\label{pri1}
\int_{\Bbb D} \left(\frac12\frac{\partial}{\partial t}|v|^2+\nu
|\o|^2\right) dx
 &=&\nu\int_{\Bbb D}  \Delta \q \, dx
  -\int_{\Bbb D }  v\cdot
 \nabla \q\,dx \n \\
 &=&\nu\int_{\partial \Bbb D }\hat{n} \cdot \nabla \q \, dS
  -\int_{\partial\Bbb D }  v\cdot
\hat{ n} \q\,dS \n \\
 &:=& I_1+I_2,
 \eq
 where $\hat{n}$ is the outward unit normal vector on $\partial \Bbb
 D$.
 \bq\label{pri2}
 I_1&=&  \nu \int_{\{\q(x,t)=\beta(t)\}}
\frac{ \nabla \q}{|\nabla \q |}\cdot \nabla  \q  dS-\nu
\int_{\{\q(x,t)=\alpha(t\})}
 \frac{ \nabla \q}{|\nabla \q |}\cdot \nabla  \q  dS\n \\
 &=& \nu \int_{\{\q(x,t)=\beta(t\})}|\nabla Q | dS-\nu \int_{\{\q(x,t)=\alpha(t\})}
|\nabla Q | dS. \eq
  Next, by the condition $\mathrm{div} \,v=0$ we obtain
\bq\label{pri3}
 I_2&=&  -\int_{\{\q(x,t)=\beta(t)\}}
 \hat{ n}\cdot v \q  dS-\int_{\{\q (x,t)=\alpha(t) \}}
  \hat{ n}\cdot v  dS\n \\
  &=&\beta(t)\int_{\{\q(x,t)=\beta(t)\}}
(-\hat{n})\cdot v  dS -\alpha(t) \int_{\{\q(x,t) =\alpha(t) \}}
  \hat{ n}\cdot v  dS\n \\
  &=&\beta(t)\int_{\{\q(x,t)>\beta(t)\}} \mathrm{div} \,v \, dx-
  \alpha(t)\int_{\{\q(x,t) >\alpha(t) \}}
 \mathrm{div} \,v \, dx=0.\n \\
 \eq
Combining (\ref{pri1})-(\ref{pri3}), we obtain (\ref{th1}).
 $\square$\\
\ \\
{\bf Proof of Theorem 1.2 } We set $\mathcal{Q}(x,t):=
Q(x,t)-p_0(t)$ here.
 Let $\alpha(t)<\beta(t)$ be  real valued functions on $(0,
T)$  with $\alpha(t) <\mathcal{ Q}(x,t)<\beta(t)$.
 The case $\alpha(t)\beta(t) >0$ is already covered  by Theorem 1.1. Here
we assume that
$$m(t):=\inf_{x\in \Bbb R^N} \mathcal{Q}(x,t)\leq \alpha(t)<0<\beta(t)\leq \sup_{x\in \Bbb R^N} \mathcal{Q}(x,t):=M(t).
$$
Since \bqn
\lefteqn{\int_{\{\alpha(t)<\mathcal{Q}(x,t)<\beta(t)\}}\left(\frac12\frac{\partial}{\partial
t}|v|^2+\nu |\o|^2\right) dx}\n \hspace{.in} \\
&&=\int_{\{\alpha(t)<\mathcal{Q}(x,t) <
M(t)\}}\left(\frac12\frac{\partial}{\partial t}|v|^2+\nu
|\o|^2\right)
dx\\
 &&\hspace{1.in}-\int_{\{\beta(t)<\mathcal{Q}(x,t)<
M(t)\}}\left(\frac12\frac{\partial}{\partial t}|v|^2+\nu
|\o|^2\right)
dx\n \\
&&= \int_{\{\alpha(t)<\mathcal{Q}(x,t) <
M(t)\}}\left(\frac12\frac{\partial}{\partial t}|v|^2+\nu
|\o|^2\right) dx + \int_{\{ \mathcal{Q}(x,t)=\beta(t)\}} |\nabla
\mathcal{Q}| dS \eqn
 by application of Theorem 1.1,
 it suffices to
show that
  \bb \label{claim}
\int_{\{\alpha(t)<\mathcal{Q}(x,t)\}}\left(\frac12\frac{\partial}{\partial
t}|v|^2+\nu |\o|^2\right) dx=-\int_{\{ \mathcal{Q}(x,t)=\alpha(t)\}}
|\nabla \mathcal{Q}| dS \ee for $m(t)<\alpha(t)<0$.
  The set $\Bbb
D^\alpha_+:=\{x\in\Bbb R^N\, |\, \alpha(t)<\mathcal{Q}(x,t) \}$ is
unbounded, while $\Bbb D^\alpha_-:=\{x\in\Bbb R^N\, |\,
\alpha(t)>\mathcal{Q}(x,t) \}$ is a bounded set. Moreover,
$$\partial
\Bbb D^\alpha_+=\partial \Bbb D^\alpha_- \cup \{ \infty\}. $$
 Let $R>0$ be large enough so that
 $ \partial
\Bbb D^\alpha_-  \subset B_R (0):=\{ x\in \Bbb R^N\,|\, |x|<R\}. $
We set $\mathcal{Q}^+_\alpha:= [\mathcal{Q} +\alpha]_+ =
\max\{\mathcal{Q} +\alpha, 0\} $ below. We introduce the radial
cut-off function $\sigma\in C_0 ^\infty(\Bbb R^N)$ such that
$$
   \sigma(|x|)=\left\{ \aligned
                  &1 \quad\mbox{if $|x|<1$},\\
                     &0 \quad\mbox{if $|x|>2$},
                      \endaligned \right.
$$
and $0\leq \sigma  (x)\leq 1$ for $1<|x|<2$.  Then, for each $R
>0$, we define
 $$
\s \left(\frac{|x|}{R}\right):=\s_R (|x|)\in C_0 ^\infty (\Bbb R^N).
$$
 Multiplying (\ref{prin}) by
$ \si $, and integrating on $\Bbb D^\alpha_+ $, we obtain by the
divergence theorem,
 \bq\label{pri13a}
 \lefteqn{\int_{\Bbb D^\alpha_+} \left(\frac12\frac{\partial}{\partial t} |v|^2+\nu |\o|^2\right)
 \si dx=\nu \int_{\Bbb D^\alpha_+} \si \Delta \mathcal{Q}^+_\alpha dx-\int_{\Bbb D^\alpha_+}
  \si v\cdot \mathcal{Q}^+_\alpha dx}\n \\
 &&=\nu\int_{\partial \Bbb D^\alpha_- } \hat{n}\cdot \nabla \mathcal{Q}^+_\alpha
\si \,
  dS-\nu \int_{\Bbb D^\alpha_+} \nabla \mathcal{Q}^+_\alpha \cdot \nabla \si
  dx
  -\int_{\partial \Bbb D^\alpha_- }  v\cdot \hat{n}\,
  \mathcal{Q}^+_\alpha\si\,dx\n \\
  &&\hspace{1.in}+\int_{ \Bbb D^\alpha_+ }\mathcal{Q}^+_\alpha v\cdot \nabla \si dx\n \\
 &&:= I_1+I_2+I_3+I_4,
 \eq
 where $\hat{n}=-\frac{\nabla \mathcal{Q}^+_\alpha}{ |\nabla \mathcal{Q}^+_\alpha|}$
 is the outward unit normal vector to $\partial \Bbb D^\alpha_- $.
As in the proof of Theorem 1.1, we have
 \bq\label{pri23}
 I_1&=& - \nu\int_{\partial \Bbb D^\alpha_-} \frac{\nabla \mathcal{Q}^+_\alpha}{ |\nabla \mathcal{Q}^+_\alpha|}\cdot
 \nabla \mathcal{Q}^+_\alpha\si  dS=- \nu\int_{\partial \Bbb D^\alpha_- } |\nabla
 \mathcal{Q}^+_\alpha | \si dS\n \\
 &=&-\nu\int_{\{\mathcal{Q}(x,t)=\alpha(t)\}} |\nabla
 \mathcal{Q}| \si dS
 =-\nu\int_{\{\mathcal{Q}(x,t)=\alpha(t)\}} |\nabla
 \mathcal{Q}|  dS,\n \\
 \eq
 since $\si=1$ on $B_R (0)\supset\partial \Bbb D^\alpha_-=\{x\in \Bbb R^N \,|\, \mathcal{Q}(x,t)=\alpha(t)\}$.
By the fact $ \mathcal{Q}^+_\alpha =0 $ on $\partial\Bbb D^\alpha_+
$ we have
 \bb\label{pri33}
 I_3=0.
 \ee
Applying the divergence theorem again, we have
  \bqn
    I_2&=&-\nu \int_{\Bbb D^\alpha_+} \nabla \mathcal{Q}\cdot \nabla \si
  dx
  =-\nu\int_{ \partial \Bbb D^\alpha_-}\mathcal{Q}\,\hat{n}\cdot \nabla \si
    dS+\nu \int_{ \Bbb D^\alpha_+}\mathcal{Q} \Delta \si dx\n \\
&=&\nu \int_{ \Bbb D^\alpha_+}\mathcal{Q} \Delta \si dx,\eqn
  since  $\nabla \si =0$ on $B_R (0)\supset\partial \Bbb D^\alpha_-$.
  Hence, we have the estimates,
  \bqn
  |I_2|&\leq& \nu\int_{\Bbb R^N} |\mathcal{Q}||\Delta \si |dx
  \leq \frac{ \|D^2 \s\|_{L^\infty}}{R^2} \|\mathcal{Q}\|_{L^{\frac{N}{N-2} }(R\leq
  | x|\leq 2R)}\left(\int_{\{R\leq
  | x|\leq 2R\}} \, dx\right)^{\frac{2}{N}}\n \\
&\leq & C\nu \|D^2\s\|_{L^\infty}
\left( \int_{\{R\leq |x|\leq 2R\}  } (|p-p_0| +\frac12 |v|^2)^{\frac{N}{N-2}}\, dx \right)^{\frac{N-2}{N}}\n \\
&\leq& C \nu \left(\|p-p_0\|_{L^{\frac{N}{N-2}} (R\leq |x|\leq 2R)}
+\|v\|_{L^{\frac{2N}{N-2}} (R\leq |x|\leq 2R)}^2\right)\to 0 \eqn as
$R\to \infty$, which follows from the assumption $v \in L^\infty(0,
T; L^{\frac{2N}{N-2}}(\Bbb R^N))$, and the Calderon-Zygmund
inequality, $ \|p-p_0\|_{L^{\frac{N}{N-2}}}\leq C
\|v\|_{L^{\frac{2N}{N-2}}}^2, $ where we used the well-known
velocity-pressure relation,
 $p-p_0=\sum_{j,k}^N R_jR_k (v_jv_k)$ with the Riesz transforms $R_j, j=1, \cdots , N$ in $\Bbb R^N$.
 For the estimate of $I_4$, since $ \partial \Bbb D_- ^\alpha \subset B_R (0)$, we have
\bq\label{j2}
|I_4|&=& \left|\int_{\{R\leq |x|\leq 2R\}  }(\mathcal{Q} +\alpha)  v\cdot \nabla \si dx\right|= \left|\int_{\{R\leq |x|\leq 2R\}  }\mathcal{Q} v\cdot \nabla \si dx\right|\n \\
&\leq & \frac{1}{R} \|D\s\|_{L^\infty}\left( \int_{\{R\leq |x|\leq 2R\}  }\, dx \right)^{\frac{1}{N}}
\left\{ \int_{\{R\leq |x|\leq 2R\}  } \left(|v||p-p_0|+ \frac12 |v|^3\right)^{\frac{N}{N-1}}\, dx \right\}^{\frac{N-1}{N}}\n \\
&\leq & C \left(\|v\|_{L^{\frac{3N}{N-1}}(R\leq |x|\leq 2R)} \|p-p_0\|_{L^{\frac{3N}{2(N-1)}}} +\|v\|_{L^{\frac{3N}{N-1}}(R\leq |x|\leq 2R)}^3\right)\n \\
&\leq& C \left(\|v\|_{L^{\frac{3N}{N-1}}(R\leq |x|\leq 2R)}
\|v\|_{L^{\frac{3N}{N-1}}} ^2 +\|v\|_{L^{\frac{3N}{N-1}}(R\leq
|x|\leq 2R)}^3\right)\to 0 \eq as $R\to \infty$ by the hypothesis
(\ref{hypo}), where we used the fact $\alpha\int_{\{R\leq |x|\leq
2R\}} v\cdot \nabla \si dx=0$ in the first line. On the other hand,
if we pass $R\to \infty$ in the left hand side of  (\ref{pri13a}),
we note that $ \int_{\Bbb D^\alpha_+}
\frac12\frac{\partial}{\partial t} |v|^2
 \si dx\to \int_{\Bbb D^\alpha_+}
\frac12\frac{\partial}{\partial t} |v|^2
 dx
$ by the dominated convergence theorem, while $ \int_{\Bbb
D^\alpha_+} |\o |^2
 \si dx\to \int_{\Bbb D^\alpha_+}
|\o |^2
 dx
$ by the monotone convergence theorem. Therefore (\ref{claim})
follows from  (\ref{pri13a}) by passing $R\to \infty$ in both sides
of the
equality. $\square$\\
 \ \\
 $$ \mbox{\bf Acknowledgements} $$
This work was supported partially by  NRF Grant no. 2006-0093854,
and by the Chung-Ang University Research Grants in 2012.

\end{document}